\documentclass[amstex,12pt, amssymb]{article}

\usepackage{mathtext}
\usepackage[cp1251]{inputenc}
\usepackage[T2A]{fontenc}
\usepackage[dvips]{graphicx}
\usepackage{amsmath}
\usepackage{amssymb}
\usepackage{amsxtra}
\usepackage{latexsym}
\usepackage{ifthen}

\newcounter{lemma}[section]

\newcounter{corollary}[section]

\newcounter{remark}[section]

\newcounter{theorem}[section]

\newcounter{proposition}[section]

\newcounter{example}

\numberwithin{equation}{section}

\begin{document}

\markboth{\centerline{E.~SEVOST'YANOV}} {\centerline{ON LOGARITHMIC
H\"{O}LDER...}}

\def\cc{\setcounter{equation}{0}
\setcounter{figure}{0}\setcounter{table}{0}}

\overfullrule=0pt


\author{{EVGENY SEVOST'YANOV}\\}

\title{
{\bf ON LOGARITHMIC H\"{O}LDER CONTINUITY OF MAPPINGS ON THE
BOUNDARY}}

\date{\today}
\maketitle

\begin{abstract}
We study mappings satisfying the so-called inverse Poletsky
inequality. Under integrability of the corresponding majorant, it is
proved that these mappings are logarithmic H\"{o}lder continuous in
the neighborhood of the boundary points. In particular, the
indicated properties hold for homeomorphisms whose inverse satisfy
the weighted Poletsky inequality.
\end{abstract}


\section{Introduction}
As is known, quasiconformal mappings and mappings with bounded
distortion are H\"{o}lder continuous with some exponent (see, e.g.,
\cite[Theorem~III.C]{Ahl}, \cite[Theorem~3.2.II]{LV},
and~\cite[Theorem~18.2, Remark~18.4]{Va}). In this regard, there are
also classical results concerning mappings with bounded distortion,
or quasiregular mappings, which are rightly called quasiconformal
mappings with branch points (see, for example,
\cite[Theorem~3.2]{MRV$_1$} and \cite[Theorem~1.1.2]{Re}). There are
many generalizations of these results to more general classes of
mappings with finite distortion. In this case, quite often, the
usual estimates of the H\"{o}lder type do not hold for mappings,
however, more general logarithmic estimates may be satisfied (see,
for example, \cite[Theorems 4 and 5]{Cr},
\cite[Theorem~7.4]{MRSY$_1$}, \cite[Theorem~3.1]{MRSY$_2$},
\cite[Theorem~5.11]{RS} and~\cite[Theorems~1.1.V and 2.1.V]{Suv}).

\medskip
Distance distortion theorems and H\"{o}lder-type estimates have been
discussed and studied in our last two articles, \cite{RSS} and
\cite{SSD}. In particular, in~\cite{SSD}, we obtained estimates for
the distortion of mappings with the inverse Poletsky inequality at
the inner points of a given domain. The main purpose of this
manuscript is to obtain similar estimates not only at the inner, but
also at the boundary points of a given domain, which, for the sake
of simplicity, will be assumed to be the unit ball.

\medskip
In what follows, for the sets $A, B\subset{\Bbb R}^n$ we set, as
usual,
$${\rm diam}\,A=\sup\limits_{x, y\in A}|x-y|\,,\quad {\rm dist}\,(A, B)=\inf\limits_{x\in A,
y\in B}|x-y|\,.$$
Let $x_0\in\overline{D},$ $x_0\ne\infty,$
$$B(x_0, r)=\{x\in {\Bbb R}^n: |x-x_0|<r\}\,, \quad {\Bbb B}^n=B(0, 1)\,.$$
Sometimes, instead of ${\rm dist}\,(A, B),$ we also write $d(A, B),$
if a misunderstanding is impossible. A Borel function $\rho:{\Bbb
R}^n\,\rightarrow [0,\infty] $ is called {\it an admissible} for a
family $\Gamma$ of paths $\gamma$ in ${\Bbb R}^n,$ if the relation
\begin{equation}\label{eq1.4}
\int\limits_{\gamma}\rho (x)\, |dx|\geqslant 1
\end{equation}
holds for any locally rectifiable path $\gamma\in\Gamma.$
{\it A modulus} of $\Gamma $ is defined as follows:
\begin{equation}\label{eq1.3gl0}
M(\Gamma)=\inf\limits_{\rho \in \,{\rm adm}\,\Gamma}
\int\limits_{{\Bbb R}^n} \rho^n (x)\,dm(x)\,.
\end{equation}
Let $Q:{\Bbb R}^n\rightarrow [0, \infty]$ be a Lebesgue measurable
function. We say that {\it $f$ satisfies the inverse Poletsky
inequality}, if the relation
\begin{equation}\label{eq2*A}
M(\Gamma)\leqslant \int\limits_{f(D)} Q(y)\cdot\rho_*^n(y)\, dm(y)
\end{equation}
holds for any family of (locally rectifiable) paths $\Gamma$ in $D$
and any $\rho_*\in {\rm adm}\,f(\Gamma).$
Note that estimates~(\ref{eq2*A}) hold in many classes of mappings
(see, e.g., \cite[Theorem~3.2]{MRV$_1$}, \cite[Theorem~6.7.II]{Ri}
and \cite[Theorem~8.5]{MRSY$_1$}). A mapping $f:D\rightarrow {\Bbb
R}^n$ is called {\it a discrete} if $\{f^{-1}\left(y\right)\}$
consists of isolated points for any $y\,\in\,{\Bbb R}^n,$ and {\it
an open,} if the image of any open set $U\subset D$ is an open set
in ${\Bbb R}^n.$ A mapping $f$ between domains $D$ and
$D^{\,\prime}$ is said to be {\it a closed} if $f(E)$ is closed in
$D^{\,\prime}$ for any closed set $E\subset D$ (see, e.g.,
\cite[section~3]{Vu$_1$}).

\medskip
Given $\delta>0,$ a non-degenerate continuum $A\subset{\Bbb B}^n$
and a Lebesgue measurable function $Q:~{\Bbb B}^n\rightarrow [0,
\infty]$ denote ${\frak S}_{\delta, A, Q }$ a family of all open,
discrete and closed mappings $f$ of the open unit ball onto itself
such that the relation~(\ref{eq2*A}) holds and, in addition, ${\rm
dist}\,(f^{\,-1}(A),
\partial {\Bbb B}^n)\geqslant~\delta.$ The following statement holds.

\medskip
\begin{theorem}\label{th1}
{\sl\,Let $Q\in L^1({\Bbb B}^n).$ Then any $f\in {\frak S}_{\delta,
A, Q }$ has a continuous extension $f:\overline{{\Bbb
B}^n}\rightarrow\overline{{\Bbb B}^n},$ and, in addition, for any
$x_0\in\partial {\Bbb B}^n$ there is $C_n>0$ and $r_0=r_0(x_0)>0$
such that
\begin{equation}\label{eq2C}
|f(x)-f(x_0)|\leqslant\frac{C_n\cdot (\Vert
Q\Vert_1)^{1/n}}{\log^{1/n}\left(1+\frac{r_0}{|x-x_0|}\right)}
\end{equation}
for $x\in B(x_0, r_0)\cap \overline{{\Bbb B}^n},$
where $\Vert Q\Vert_1$ denotes the $L^1$-norm of the function $Q$ in
$L^1({\Bbb B}^n).$
 }
\end{theorem}

\medskip
\section{Auxiliary lemma and proof of Theorem~\ref{th1}} Before
by proving the basic statement we prove the following important
lemma.

\medskip
\begin{lemma}\label{lem1}
{\sl\, Let $E$ be a continuum in~${\Bbb B}^n.$ Now, there is
$\delta_1>0$ such that ${\frak S}_{\delta, A, Q }\subset {\frak
S}_{\delta_1, E, Q }.$ In other words, if $f$ is an open discrete
and closed mapping of the unit ball onto itself with condition
(\ref{eq2*A}), such that ${\rm dist}\,(f^{\,-1}(A),\partial{\Bbb
B}^n)\geqslant~\delta,$ then there is $\delta_1> 0,$ independent on
$f$ such that ${\rm dist}\,(f^{\,-1}(E),\partial{\Bbb B}^n)\geqslant
~\delta_1.$
 }
\end{lemma}

\medskip
\begin{proof}
Let us prove this statement by contradiction. Suppose that the
conclusion of the lemma is not correct. Then there are sequences
$y_m\in E,$ $f_m\subset{\frak S}_{\delta, A, Q }$ and $x_m\in {\Bbb
B}^n$ such that $f_m(x_m)=y_m$ and ${\rm dist}\,(x_m, \partial {\Bbb
B}^n)\rightarrow 0$ as $m\rightarrow\infty.$ Without loss of
generality, we may assume that $x_m\rightarrow x_0$ as
$m\rightarrow\infty.$ By~\cite[Theorem~3.1]{SSD} $f_m$ has a
continuous extension at $x_0,$ moreover, $\{f_m\}_{m=1}^{\infty}$ is
equicontinuous at~$x_0$ (see, e.g., \cite[Theorem~1.2]{SSD}). Now,
for any $\varepsilon>0$ there is $m_0\in {\Bbb N}$ such that
$|f_m(x_m)-f_m(x_0)|<\varepsilon$ for $m\geqslant m_0.$ On the other
hand, $f_m(x_0)\in {\Bbb S}^{n-1}$ because $f_m$ is closed. Thus, by
the triangle inequality, $|f_m(x_m)-f_m(x_0)|\geqslant
1-|f_m(x_m)|\geqslant 1-\delta_0,$ where $\delta_0=\sup\limits_{x\in
E}|x|.$ Finally, we have a contradiction, because
$|f_m(x_m)-f_m(x_0)|\geqslant 1-\delta_0$ and, simultaneously
$|f_m(x_m)-f_m(x_0)|<\varepsilon$ as $m\geqslant m_0.$ The resulting
contradiction refutes the assumption made above. Lemma is
proved.~$\Box$
\end{proof}

\medskip
{\it Proof of Theorem~\ref{th1}}. The possibility of a continuous
extension of $f$ to the boundary of ${\Bbb B}^n$ is established
in~\cite[Theorem~3.1]{SSD}. In particular, the weakly flatness of
$\partial {\Bbb B}^n={\Bbb S}^{n-1}$ follows by~\cite[Theorems~17.10
and 17.12]{Va}.

\medskip
Let us prove the logarithmic H\"{o}lder continuity~(\ref{eq2C}). It
suffices to prove relation~(\ref{eq2C}) for the case when $x\in
{\Bbb B}^n\cap B(x_0, r_0),$ since the general case $x\in
\overline{{\Bbb B}^n}\cap B(x_0, r_0)$ is attained by passing to the
limit as $x\rightarrow x_*,$ where $x_*\in \partial {\Bbb B}^n\cap
B(x_0, r_0).$ Let~$x_0 \in
\partial{\Bbb B}^n,$ let $0<\delta<1$ and let $E=B(0,
\delta/2)\subset {\Bbb B}^n.$ By Lemma~\ref{lem1}, there
is~$\delta_1>0$ such that~${\rm dist}\,(f^{\,-1}(E),
\partial {\Bbb B}^n)\geqslant \delta_1$ for any $f\in {\frak S}_{\delta, A, Q}.$
By~\cite[Theorem ~1.2]{SSD}, the family~${\frak S}_{\delta, A, Q}$
is equicontinuous in~$\overline{{\Bbb B}^n}.$ Thus, for a number
$0<\delta_0<1/4$ there is a neighborhood $U\subset B(x_0,
\delta_1/2)$ of $x_0$ such that $|f(x)-f(y)|<\delta_0$ for any $x,
y\in U\cap {\Bbb B}^n$ and $f\in{\frak S}_{\delta, A, Q}.$ Let
$f(x)\ne f(y)$ and
$$\varepsilon_0:=|f(x)-f(y)|<\delta_0\,.$$
Let us join the points $f(x)$
and $0$ by segment $I.$ The points  $f(x),$ $0$ and $f(y)$ form the
plane $P.$ Consider a circle
$$S=\{z\in P:
|z-f(x)|=\varepsilon_0\}\,.$$
The position of the point $z=f(y)$ on the circle $S$ is completely
determined by the angle $\varphi,$ $\varphi\in [-\pi, \pi),$ between
the vector $-f(x)$ and the radius-vector of the point $z.$ The
points on the circle are denoted further in the polar coordinates by
the pairs $z=(\varepsilon_0, \varphi).$ Our further goal is to
investigate the main three cases regarding the intervals of change
of this angle.

\begin{figure}[h]
\centerline{\includegraphics[scale=0.4]{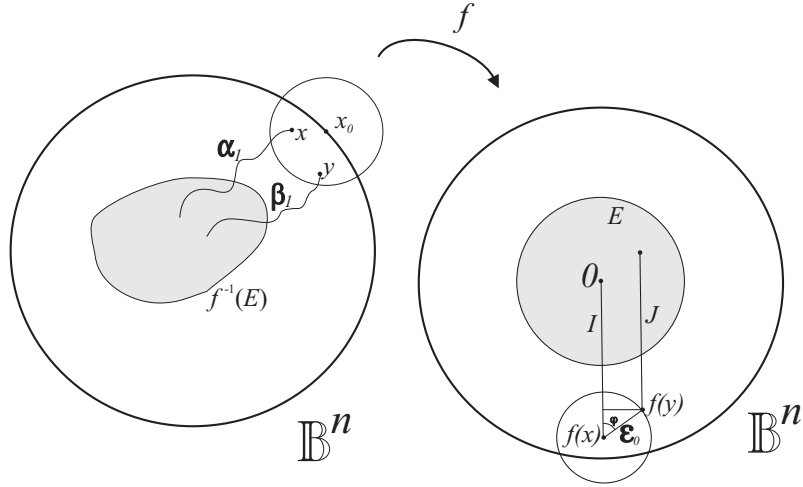}} \caption{To
the proof of Theorem~\ref{th1}, the case~{\textbf{1}}}\label{fig1}
\end{figure}

{\bf Case 1.} ''Large angles'':  $\varphi\in [\pi/4, 3\pi/4],$ or
$\varphi\in [-\pi/4, -3\pi/4].$ Let $\varphi\in [\pi/4, 3\pi/4].$
Consider the ray
$$r=r(t)=f(y)+te\,,\quad
e=-f(x)/|f(x)|\,,\quad t>0.\,$$
By construction, the ray $r$ is parallel to the segment $I.$ For
$t=|f(x)|,$ we have $r(|f(x)|)=f(y)-f(x)$ and
$|r(|f(x)|)|=\varepsilon_0<\delta_0,$ i.e., the point $r(|f(x)|)$
belongs to $E.$ Let $J$ be a segment of the ray $r,$ contained
between the points $f(y)$ and $r(|f(x)|).$ The distance between $I$
and $ J $ is calculated as follows:
\begin{equation}\label{eq1}
{\rm dist}\,(I, J)=\varepsilon_0\sin\varphi\geqslant
\frac{\sqrt{2}}{2}\varepsilon_0\,,
\end{equation}
see Figure~\ref{fig1} for illustration.
Similarly, the situation $\varphi\in [-\pi/4, -3\pi/4]$ is
considered. We may show that in this case the formula~(\ref{eq1})
also holds.

\medskip
Let $\alpha_1$ and $\beta_1$ be (total) $f$-liftings of paths $I$
and $J$ starting at points $x$ and $y,$ correspondingly. These
liftings are well defined and exist by~\cite[Lemma~3.7]{Vu$_1$}.
Now, by the definition, $|\alpha_1|\cap f^{\,-1}(E)\ne\varnothing\ne
|\beta_1|\cap f^{\,-1}(E).$ Since ${\rm dist}\,(f^{\,-1}(E),
\partial {\Bbb B}^n)\geqslant \delta_1$ and $x, y\in B(x_0, \delta_1/2),$
\begin{equation}\label{eq4}
{\rm diam\,}(\alpha_1)\geqslant \delta_1/2\,,\quad  {\rm
diam\,}(\beta_1)\geqslant \delta_1/2\,.
\end{equation}
Let
$$\Gamma:=\Gamma(\alpha_1, \beta_1, {\Bbb B}^n)\,.$$
Then on the one hand, by~\cite[Lemma~4.3]{Vu$_2$}
\begin{equation}\label{eq7A}
M(\Gamma)\geqslant (1/2)\cdot M(\Gamma(\alpha_1, \beta_1, {\Bbb
R}^n))\,,
\end{equation}
and on the other hand, by~\cite[Lemma~7.38]{Vu$_3$}
\begin{equation}\label{eq7B}
M(\Gamma(\alpha_1, \beta_1, {\Bbb R}^n))\geqslant
c_n\cdot\log\left(1+\frac1m\right)\,,
\end{equation}
where $c_n>0$ is come constant depending only on $n,$
$$m=\frac{{\rm dist}(\alpha_1, \beta_1)}{\min\{{\rm diam\,}(\alpha_1),
{\rm diam\,}(\beta_1)\}}\,.$$
Combining~(\ref{eq4}) and~(\ref{eq7B}), and taking into account that
${\rm dist}\,(\alpha_1, \beta_1)\leqslant |x-y|,$ we obtain that
\begin{equation}\label{eq7C}
M(\Gamma)\geqslant \widetilde{c_n}\cdot
\log\left(1+\frac{\delta_1}{2{\rm dist}(\alpha_1,
\beta_1)}\right)\geqslant \widetilde{c_n}\cdot
\log\left(1+\frac{\delta_1}{2|x-y|}\right)\,,
\end{equation}
where $\widetilde{c_n}>0$ is some constant depending only on~$n.$

\medskip
We now obtain an upper bound for $M(\Gamma).$ Set
$$\rho(x)= \left\{
\begin{array}{rr}
\frac{\sqrt{2}}{\varepsilon_0}, & x\in {\Bbb B}^n,\\
0,  &  x\not\in {\Bbb B}^n\,.
\end{array}
\right. $$
By~(\ref{eq1}) $\rho$ satisfies the relation~(\ref{eq1.4}) for the
family $f(\Gamma).$ Then by the definition of the family ${\frak
S}_{\delta, A, Q }$ we obtain that
\begin{equation}\label{eq14***}
M(\Gamma)\leqslant \frac{2^{n/2}}{\varepsilon_0^n}\int\limits_{{\Bbb
B}^n} Q(y)\,dm(y)=2^{n/2}\cdot \frac{\Vert
Q\Vert_1}{{|f(x)-f(y)|}^n}\,,
\end{equation}
where $\Vert Q\Vert_1$ denotes the $L^1$-norm of the function $Q$ in
${\Bbb B}^n.$ By~(\ref{eq7C}) and (\ref{eq14***}) we obtain that
$$\widetilde{c_n}\cdot \log\left(1+\frac{\delta_1}{2|x-y|}\right)\leqslant
2^{n/2}\cdot\frac{\Vert Q\Vert_1}{{|f(x)-f(y)|}^n}\,.$$
The desired inequality~(\ref{eq2C}) follows from the last relation
by passing to the limit as $y\rightarrow x_0\,,$ where
$C_n:=2^{n/2}\cdot\widetilde{c_n}^{-1/n}$ and $r_0=\min\{\delta_1/2,
d(x_0, \partial U)\}.$

\medskip
{\bf Case~2.} \begin{figure}[h]
\centerline{\includegraphics[scale=0.4]{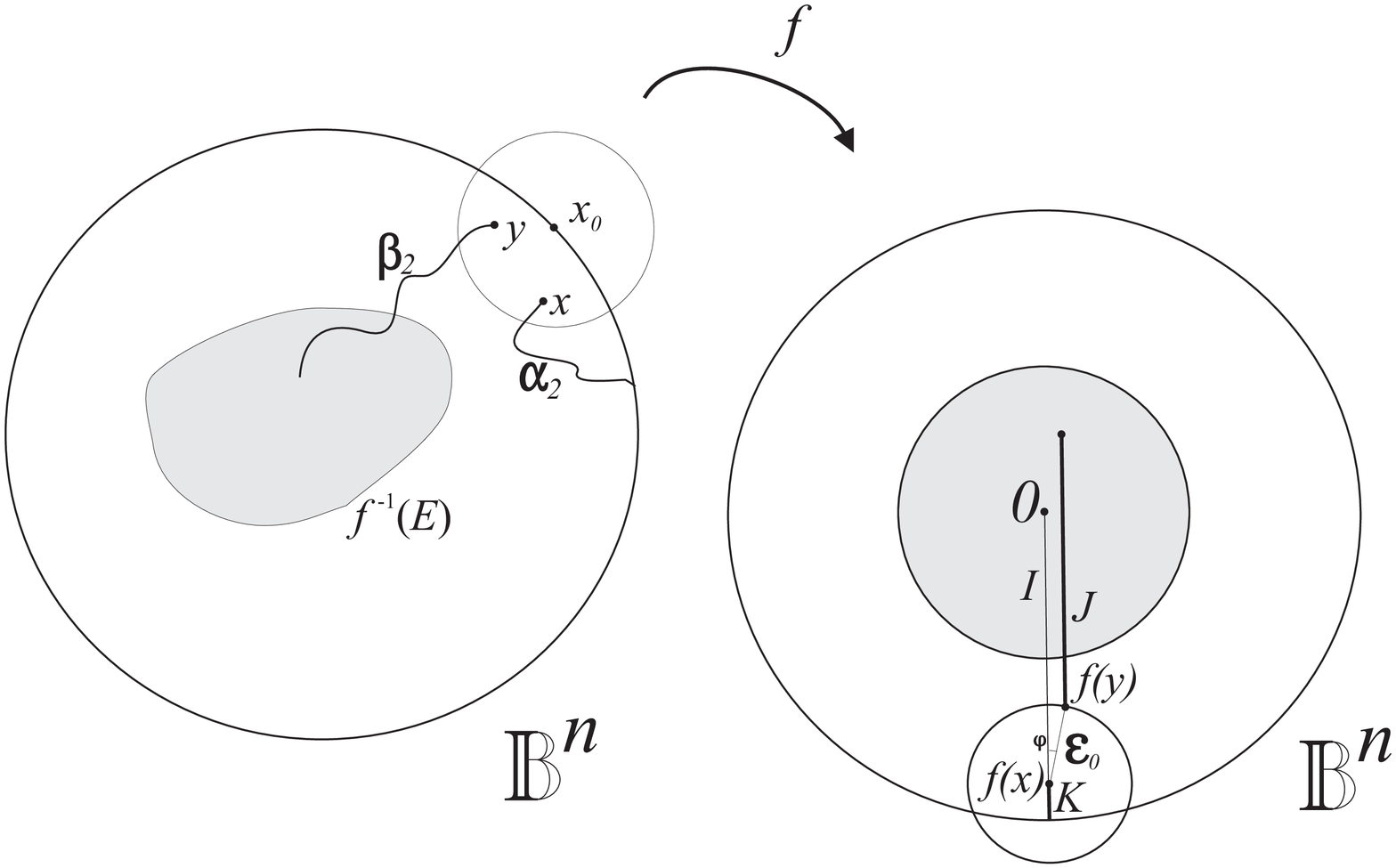}} \caption{To
the proof of Theorem~\ref{th1}, the case~{\textbf{2}}}\label{fig2}
\end{figure}
''Small angles'': $\varphi\in [-\pi/4, \pi/4),$ see
Figure~\ref{fig2} for the illustration. In this case, we denote by
$K$ the segment joining the point $f(x)$ with the unit sphere in the
direction of the vector $f(x)/|f(x)|.$ Then
\begin{equation}\label{eq2}
{\rm dist}\,(K, J)=\varepsilon_0\,.
\end{equation}
Let $\alpha_2$ be a maximal $f$-lifting of $K^{\,\prime}$ starting
at $x,$ where $K^{\,\prime}$ is obtained from $K$ by discarding its
endpoint lying on the unit sphere. Such a lift exists and tends with
its end to~$\partial {\Bbb B}^n={\Bbb S}^{n-1}$ (see, e.g.,
\cite[Lemma~3.12]{MRV$_2$}). Let also $\beta_2$ be a total
$f$-lifting of $J$ starting at the point $y$ (such a lifting exists
by \cite[Lemma~3.7]{Vu$_1$}). Arguing similarly to a case~1, we
obtain that
\begin{equation}\label{eq4A*}
{\rm diam\,}(\beta_2)\geqslant \delta_1/2\,.
\end{equation}
Let
$$\Gamma:=\Gamma(\alpha_2\cup {\Bbb S}^{n-1}, \beta_2, \overline{{\Bbb B}^n})\,.$$
According to the above, $\overline{\alpha_2}\setminus
\alpha_2\subset {\Bbb S}^{n-1}.$ Thus $\overline{\alpha_2\cup {\Bbb
S}^{n-1}}=\overline{\alpha_2}\cup {\Bbb S}^{n-1}=\alpha_2\cup {\Bbb
S}^{n-1},$ consequently, $\alpha_2\cup {\Bbb S}^{n-1}$ is a closed
set. Observe that $\overline{\alpha_2}$ is connected (see, e.g.,
\cite[Corollary~3(ii).II.46.5]{Ku}). Moreover, the set $\alpha_2\cup
{\Bbb S}^{n-1}=\overline{\alpha_2}\cup {\Bbb S}^{n-1}$ is connected
as the union of two connected sets $\overline{\alpha_2}$ and ${\Bbb
S}^{n-1},$ that have at least one common point (see, e.g.,
\cite[Corollary~3(i).II.46.5]{Ku}).

\medskip
In this case, on the one hand, by~\cite[Lemma~4.2]{Vu$_2$}
\begin{equation}\label{eq7A*}
M(\Gamma)\geqslant (1/2)\cdot M(\Gamma(\alpha_2\cup {\Bbb S}^{n-1},
\beta_2, {\Bbb R}^n))\,,
\end{equation}
and on the other hand,  by~\cite[Lemma~7.38]{Vu$_3$}
\begin{equation}\label{eq7B*}
M(\Gamma(\alpha_2\cup {\Bbb S}^{n-1}, \beta_2, {\Bbb R}^n))\geqslant
c_n\cdot\log\left(1+\frac1m\right)\,,
\end{equation}
where $c_n>0$ is some constant depending only on~$n,$
$$m=\frac{{\rm dist}(\alpha_2\cup {\Bbb S}^{n-1}, \beta_2)}
{\min\{{\rm diam\,}(\alpha_2\cup {\Bbb S}^{n-1}), {\rm
diam\,}(\beta_2)\}}\,.$$
Then combining~(\ref{eq4A*}) and~(\ref{eq7B*}), and taking into
account that ${\rm dist}\,(\alpha_2\cup {\Bbb S}^{n-1},
\beta_2)\leqslant |x-y|,$ we obtain that
\begin{equation}\label{eq7C*}
M(\Gamma)\geqslant \widetilde{c_n}\cdot
\log\left(1+\frac{\delta_1}{2{\rm dist}(\alpha_2\cup {\Bbb S}^{n-1},
\beta_2)}\right)\geqslant \widetilde{c_n}\cdot
\log\left(1+\frac{\delta_1}{2|x-y|}\right)\,,
\end{equation}
where $c_n>0$ is some constant depending only on~$n.$

\medskip
Let us now establish an upper bound for $M(\Gamma).$ First of all,
note that
$$\Gamma(\alpha_2\cup {\Bbb S}^{n-1}, \beta_2, \overline{{\Bbb
B}^n})\supset \Gamma(\alpha_2\cup {\Bbb S}^{n-1}, \beta_2, {\Bbb
B}^n)$$ and $$\Gamma(\alpha_2\cup {\Bbb S}^{n-1}, \beta_2,
\overline{{\Bbb B}^n})>\Gamma(\alpha_2\cup {\Bbb S}^{n-1}, \beta_2,
{\Bbb B}^n)\,.$$ Therefore, by the principle of minority and in view
of the monotonicity of the module we obtain that
\begin{equation}\label{eq8}
M(\Gamma(\alpha_2\cup {\Bbb S}^{n-1}, \beta_2, \overline{{\Bbb
B}^n}))=M(\Gamma(\alpha_2\cup {\Bbb S}^{n-1}, \beta_2, {\Bbb
B}^n))=M(\Gamma)\,.
\end{equation}
Set
$$\rho(x)= \left\{
\begin{array}{rr}
\frac{1}{\varepsilon_0}, & x\in {\Bbb B}^n,\\
0,  &  x\not\in {\Bbb B}^n\,.
\end{array}
\right. $$
Observe that, by~(\ref{eq2}), $\rho$ satisfies the
relation~(\ref{eq1.4}) for $f(\Gamma(\alpha_2\cup {\Bbb S}^{n-1},
\beta_2, {\Bbb B}^n)).$ In this case, by~(\ref{eq8}) and the
definition of the family ${\frak S}_{\delta, A, Q }$ we obtain that
\begin{equation}\label{eq14***A}
M(\Gamma)=M(\Gamma(\alpha_2\cup {\Bbb S}^{n-1}, \beta_2, {\Bbb
B}^n))\leqslant \frac{1}{\varepsilon_0^n}\int\limits_{{\Bbb B}^n}
Q(y)\,dm(y)=\frac{\Vert Q\Vert_1}{{|f(x)-f(y)|}^n}\,.
\end{equation}
By~(\ref{eq7C*}) and~(\ref{eq14***A}) we obtain that
$$\widetilde{c_n}\cdot \log\left(1+\frac{\delta_1}{2|x-y|}\right)\leqslant
\frac{\Vert Q\Vert_1}{{|f(x)-f(y)|}^n}\,.$$
The desired inequality~(\ref{eq2C}) follows from the last relation,
which is achieved by passing to the limit as~$y\rightarrow x_0\,,$
where $C_n:=\widetilde{c_n}^{-1/n}$ and $r_0=\min\{\delta_1/2,
d(x_0, \partial U)\}.$

\medskip
{\bf Case~3.}  \begin{figure}[h]
\centerline{\includegraphics[scale=0.4]{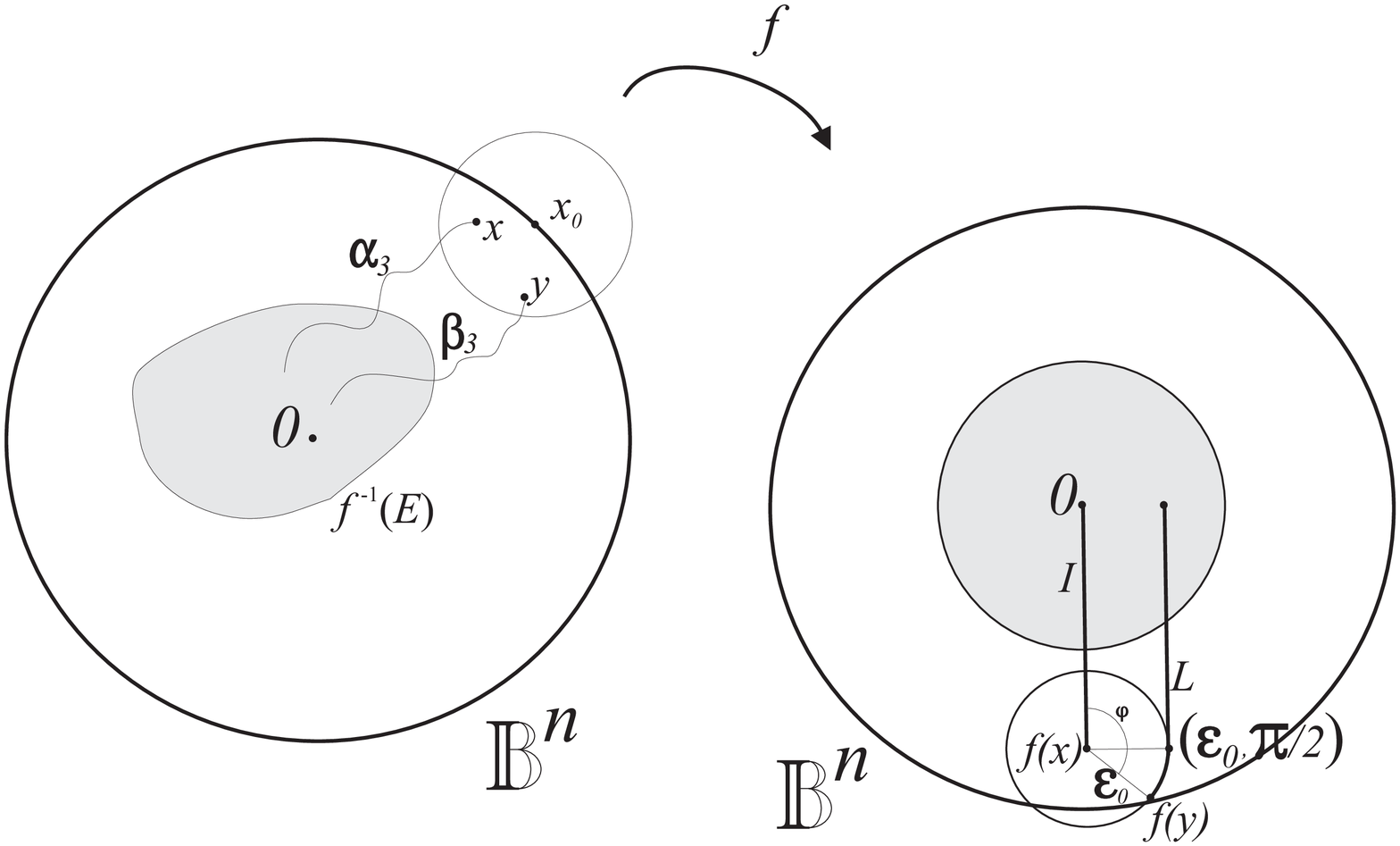}} \caption{ To
the proof of Theorem~\ref{th1}, the case~{\textbf{3}}}\label{fig4}
\end{figure} <<Very large angles>>: either $\varphi\in (3\pi/4,
\pi],$ or $\varphi\in (-\pi, -3\pi/4),$ see Figure~\ref{fig4} for
the illustration. Let, for example, $\varphi\in (3\pi/4, \pi].$
Using the methods of elementary geometry, one can show that points
$z=(\varepsilon_0, \varphi),$ having the largest absolute values of
the number $\varphi,$ also have the largest Euclidean modulus of the
number $z$. Thus, if $z=(\varepsilon_0, \varphi)\in {\Bbb B}^n,$
then the whole arc of the circle
$$z=z(\varphi)=\{(\varepsilon_0, \theta): \theta\in [\pi/2,
\varphi]\}$$ also belongs to ${\Bbb B}^n.$ Let us denote by $L_1$
the line segment
$$l=l(t)=(\varepsilon_0, \pi/2)+te,\quad
e=-f(x)/|f(x)|\,,$$
corresponding to the values of the parameter $t\in [0, |f(x)|].$
Now, we obtain that $l(|f(x)|)=(\varepsilon_0, \pi/2)-f(x),$ while
$|l(|f(x)|)|=\varepsilon_0$ by the construction. Then $L_1$
intersects the ball $E.$ Set $L=z(\varphi)\cup L_1.$ We obtain that
\begin{equation}\label{eq3}
{\rm dist}\,(L, I)=\varepsilon_0\,.
\end{equation}
It can be shown that relation~(\ref{eq3}) is also satisfied
for~$\varphi\in (-\pi, -3\pi/4).$
Let $\alpha_3$ and $\beta_3$ be $f$-liftings of $I$ and $L$ starting
$x$ and $y,$ correspondingly. As before, the existence of such
liftings is due to~\cite[Lemma~3.7]{Vu$_1$}. Set
$\Gamma=\Gamma(\alpha_3, \beta_3, {\Bbb B}^n).$ Further reasoning is
similar to what was carried out in case~1, namely, the existence of
continua $I$ and $L$ with distance $\varepsilon_0$ implies that
$$M(\Gamma)\geqslant \widetilde{c_n}\cdot
\log\left(1+\frac{\delta_1}{2{\rm dist}(\alpha_3,
\beta_3)}\right)\geqslant \widetilde{c_n}\cdot
\log\left(1+\frac{\delta_1}{2|x-y|}\right)$$
and $$M(\Gamma)\leqslant \frac{1}{\varepsilon_0^n}\int\limits_{{\Bbb
B}^n} Q(y)\,dm(y)=\frac{\Vert Q\Vert_1}{{|f(x)-f(y)|}^n}\,.$$
Hence we obtain the inequality
$$\widetilde{c_n}\cdot \log\left(1+\frac{\delta_1}{2|x-y|}\right)\leqslant
\frac{\Vert Q\Vert_1}{{|f(x)-f(y)|}^n}\,.$$
From the last inequality we obtain the statement of the theorem due
to the passage to the limit as $y\rightarrow x_0.$~$\Box$

\medskip
The analog of Theorem~\ref{th1} is also valid for mappings with a
fixed point of the unit ball. In order to formulate and prove the
corresponding statement, we introduce the following definition. For
elements $a, b\in {\Bbb B}^n$ and the Lebesgue measurable function
$Q:{\Bbb B}^n\rightarrow [0, \infty]$ denote by ${\frak F}_{a, b,
Q}$ the family of all open discrete and closed mappings $f$ of the
unit sphere onto itself, such that~$f(a)=b.$ The following statement
is true.

\medskip
\begin{theorem}\label{th2}
{\sl\,Let $Q\in L^1({\Bbb B}^n).$ Then any mapping $f\in {\frak
S}_{\delta, A, Q }$ has a continuous extension $f:\overline{{\Bbb
B}^n}\rightarrow\overline{{\Bbb B}^n},$ in addition, for any
$x_0\in\partial {\Bbb B}^n$ there exist $C_n>0$ and $r_0=r_0(x_0)>0$
such that the relation~(\ref{eq2C}) holds. }
\end{theorem}

\medskip
\begin{proof}
The possibility of a continuous extension of the mapping $f$ to
${\Bbb S}^{n-1}$ follows from Theorem~3.1 in \cite {SSD}. We prove
the logarithmic H\"{o}lder continuity of the family of extended
mappings. Let $E=\overline{B(0, 1/2)}.$ The following two cases are
possible:

\medskip
1) there exists $\delta>0$ such that ${\rm dist}\,(f^{\,-1}(E),
\partial {\Bbb B}^n)\geqslant\delta $ for any $f\in {\frak
S}_{\delta, A, Q }.$ In this case, the desired statement follows
from Theorem~\ref{th1};

\medskip
2) there exist $f_m\in {\frak S}_{\delta, A, Q }$ and $x_m, y_m\in
{\Bbb B}^n,$ $m=1,2,\ldots ,$ such that $f_m(x_m)=y_m,$ $y_m\in E$
and ${\rm dist\,}(x_m, {\Bbb S}^{n-1})\rightarrow 0$ as
$m\rightarrow\infty.$ Then, arguing exactly as in the proof of
Lemma~\ref{lem1}, we come to the conclusion that the family of
mappings~${\frak S}_{\delta, A, Q }$ is not equicontinuous at least
at one point~$x_0\in {\Bbb S}^{n-1},$ which contradicts the
assertion of Theorem~7.1 in~\cite{SSD}.

\medskip
Thus, a case 2) is impossible, and a case 1) gives the desired
statement of the theorem.~$\Box$
\end{proof}

\medskip
In addition to our article, we will illustrate our results with some
simple examples.

\medskip
\begin{example}\label{ex1}
Consider a family of plane mappings $f_n(z)=z^n,$ $n=1,2,\ldots \,,$
$z\in {\Bbb B}^2=\{z\in {\Bbb C}: |z|<1\}.$ The mappings $f_n$ have
a bounded distortion as smooth mappings whose dilatation is equal to
one. So $f_n$ satisfy inequality~(\ref{eq2*A}) for $Q(z)=N(f_n,
{\Bbb B}^2),$ where, as usual, $N$ is a multiplicity function
determined by the ratios
$$N(y, f, {\Bbb B}^2)\,=\,{\rm
card}\,\left\{z\in {\Bbb B}^2: f(z)=y\right\}\,, \qquad N(f, {\Bbb
B}^2)\,=\,\sup\limits_{y\in{\Bbb C}}\,N(y, f, {\Bbb B}^2)$$
(see, e.g., \cite[Theorem~3.2]{MRV$_1$} or
\cite[Theorem~6.7.II]{Ri}). All mappings $f_n$ are discrete and
open, in addition, preserve the boundary of the unit disk and,
therefore, are closed (see, e.g., \cite[Theorem~3.3]{Vu$_1$}). The
mappings $f_n$ also fix the point $0,$ so they satisfy all the
conditions of Theorem~\ref{th2} except one: the existence of the
integrable function $Q$ in~(\ref{eq2*A}) independent on $n.$ As a
result, the family $f_n$ is not equicontinuous at the boundary of
the unit disk.
\end{example}

\medskip
\begin{example}\label{ex2}
A simple example of conformal automorphisms
$f_n(z)=\frac{z-\frac{n-1}{n}}{1-z\frac{n-1}{n}}$ of the unit disk
on itself, for which $Q(z)\equiv 1 $ (see, e.g.,
\cite[Theorem~3.2]{MRV$_1$} or \cite[Theorem~6.7.II]{Ri}), shows
that the violation of the condition  ${\rm dist}\,(f^{\,-1}(A),
\partial {\Bbb B}^n)\geqslant~\delta$ in Theorem~\ref{th1} is
also an obstacle to performing the desired inequality~(\ref{eq2C})
on the boundary of the unit disk. We observe that the family~$f_n$
is not even equicontinuous, and thus more, logarithmically
H\"{o}lder continuous at the boundary points of the unit disk.
\end{example}

\medskip
\medskip
{\bf \noindent Evgeny Sevost'yanov} \\
{\bf 1.} Zhytomyr Ivan Franko State University,  \\
40 Bol'shaya Berdichevskaya Str., 10 008  Zhytomyr, UKRAINE \\
{\bf 2.} Institute of Applied Mathematics and Mechanics\\
of NAS of Ukraine, \\
1 Dobrovol'skogo Str., 84 100 Slavyansk,  UKRAINE\\
esevostyanov2009@gmail.com

\end{document}